\documentclass{ifacconf}

\usepackage{graphicx}      
\usepackage{natbib}        

\usepackage{amsmath}
\usepackage{amssymb}

\usepackage{algorithm}
\usepackage{savesym}
\savesymbol{AND}
\savesymbol{OR}
\savesymbol{NOT}
\savesymbol{TO}
\savesymbol{COMMENT}
\savesymbol{BODY}
\savesymbol{IF}
\savesymbol{ELSE}
\savesymbol{ELSIF}
\savesymbol{FOR}
\savesymbol{WHILE}
\usepackage{algorithmic}
\usepackage{xcolor}

\usepackage{caption}
\usepackage{subcaption}
\allowdisplaybreaks

\allowdisplaybreaks

\begin{document}
\begin{frontmatter}

\title{Distributed Generalized Wirtinger Flow for Interferometric Imaging on Networks\thanksref{footnoteinfo}} 

\thanks[footnoteinfo]{This work was partially supported by the NSF Grant CNS-1956297, the ARPA-H Strategic Initiative Seed Fund \#916012, and the Sustainable Futures Fund \#919027.}

\author{Sean M. Farrell  } 
\author{Ashok Veeraraghavan \quad } 
\author{Ashutosh Sabharwal }
\author{C\'esar A. Uribe} 

\address{Department
of Electrical and Computer Engineering, Rice University, Houston,
TX 77005, USA \\
(e-mail: \{smf5,~vashok,~ashu,~cauribe\}@rice.edu).}

\begin{abstract}                
We study the problem of decentralized interferometric imaging over networks, where agents have access to a subset of local radar measurements and can compute pair-wise correlations with their neighbors. We propose a primal-dual distributed algorithm named \emph{Distributed Generalized Wirtinger Flow} (DGWF). We use the theory of low rank matrix recovery to show when the interferometric imaging problem satisfies the Regularity Condition, which implies the Polyak-\L{}ojasiewicz inequality. Moreover, we show that DGWF converges geometrically for smooth functions. Numerical simulations for single-scattering radar interferometric imaging demonstrate that DGWF can achieve the same mean-squared error image reconstruction quality as its centralized counterpart for various network connectivity and size.
\end{abstract}


\end{frontmatter}

\section{Introduction}

\begin{figure}[ht]
    \centering
    \begin{minipage}[b]{0.55\columnwidth}
    \centering
       \subfloat[]{\includegraphics[width=1.0\columnwidth]{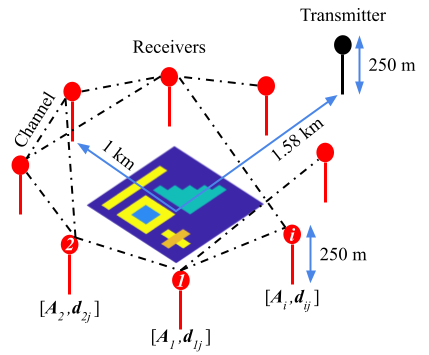}}
    \end{minipage}\hfill
     \begin{minipage}[b]{0.45\columnwidth}
     \centering
       \subfloat[]{\includegraphics[width=1.0\columnwidth]{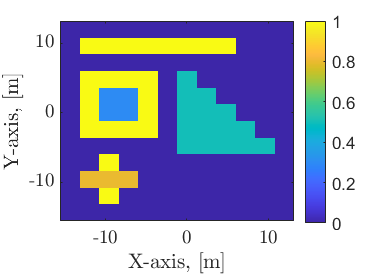}}
    \end{minipage}\hfill
    \caption{Illustration of multistatic radar interferometric imaging scene (a). In the distributed network each agent $i$ can access a subset of the data (i.e. $[\boldsymbol{A}_i,\boldsymbol{d}_{ij}]$). Scene reflectivity function used in numerical simulations (b), $12 \times 12$ image with a pixel spacing of 2.4 m.}
    \label{fig:intro}
\end{figure}

The interferometric imaging problem consists of finding an unknown complex signal from cross-correlation measurements. This nonconvex problem arises in many applications such as radar and sonar interferometry, passive electromagnetic imaging, seismic imaging, and radio astronomy~\citep{yazici2019,yazici2020,nanzer2019,papanicolaou2009}. As the scale of these applications increases, solving the interferometric imaging problem requires distributed or decentralized optimization algorithms. In practice, solving the interferometric imaging problem in a distributed manner may be beneficial when 1) sensing devices have limited memory, 2) there are data privacy concerns, 3) communication to a centralized node is unfeasible due to communication constraints. Figure \ref{fig:intro} shows a distributed radar interferometric imaging scenario, where a single transmitter illuminates a scene with a static reflectivity function. The edges (dashed lines) represent the communication channels between receivers. The $i^{th}$ receiver can only access its sampling matrix $\boldsymbol{A}_{i}$ and compute cross-correlation measurements $\boldsymbol{d}_{ij}$ with the $j^{th}$ receiver if a communication channel exits. 


Many algorithms used to solve the centralized interferometric imaging problem take inspiration from phase retrieval problems~\citep{candes2011,candes2015,duchi2018}. A gradient-descent-based iterative low rank matrix recovery approach was proposed in~\cite{yazici2015}; which \emph{lifts} the solution set to convexify the problem, leading to increased complexity. Recently, \cite{yazici2019} and \cite{yazici2020} proposed a generalized Wirtinger flow (GWF) algorithm to solve the interferometric problem, which operates in the signal domain resulting in improved computation and memory efficiencies. They show that the GWF algorithm can achieve a linear convergence rate under specific conditions on the measurements~\citep{yazici2019}. 
\cite{demanet2017} solve the centralized interferometric imaging problem using a lifting approach with a graph to encode the available interferometric measurements. They establish a relationship between graph connectivity and the robustness of recovery. 

Recently, distributed optimization techniques have been proposed to solve the phase retrieval problem. \cite{hong2018} introduced the distributed Wirtinger flow algorithm (DWF) to solve the phase retrieval problem using a perturbed proximal primal-dual approach. The authors in \cite{hong2018} show that DWF can converge to an approximate solution at a sublinear rate. \cite{garcia2021} introduced a distributed subgradient method to solve weakly convex and non-smooth problems such as the phase retrieval problem. They report a linear convergence rate, using the gradient norm as a stopping criterion when the function is locally sharp near a minimizer. 

In this work, to the best of our knowledge \citep{demanet2017,yazici2019,yazici2020}, we present the first distributed nonconvex primal-dual optimization algorithm to solve the interferometric imaging problem. \textit{The contributions of this paper are}:

\begin{itemize}
    \item We show that the Polyak-\L{}ojasiewicz (PL) inequality is satisfied for the interferometric imaging problem when its lifted forward model meets the restricted isometry properties of rank-1 positive semi-definite matrices with a sufficiently small restricted isometry constant.
    \item We show that after proper initialization the Distributed Generalized Wirtinger Flow (DGWF) algorithm converges linearly to a global optimum when the cost function is smooth and PL inequality is satisfied.
    \item We provide numerical evidence that the DGWF algorithm can achieve comparable image reconstruction quality to centralized methods \citep{yazici2019,yazici2020} for radar interferometric imaging applications.
\end{itemize}
This paper is organized as follows. Section \ref{sec:problem_form} provides the problem formulation. Section \ref{sec:DGWF} introduces the proposed DGWF algorithm. Section \ref{sec:conv_analysis} shows linear convergence for the DGWF algorithm. Numerical simulations are in Section \ref{sec:sims}, and Section \ref{sec:conclusion} concludes this paper.

\section{Problem Formulation}\label{sec:problem_form}

Suppose the desired signal is denoted by $x \in \mathbb{C}^K$, the interferometric imaging problem can be formulated as:
\begin{align}\label{eq:2.1}
    &\textrm{find} \quad\quad x \notag\\
    &\textrm{subject to} \quad d_{ij}^s = \langle a_i^s,x \rangle \overline{\langle a_j^s,x \rangle}, \notag\\
    &\textrm{for}\quad i=1,2,...,N;  j \neq i;s = 1,2,...,S  \notag
\end{align}
where $d_{ij}^s \in \mathbb{C}$ is the $s^{th}$ cross-correlation sample for the $i^{th}$ and $j^{th}$ sensing processes, $a_i^s \in \mathbb{C}^K$ is the $s^{th}$ sampling vector for the $i^{th}$ sensing process, $\overline{(\cdot)}$ is the complex conjugate, and $\langle .,. \rangle$ denotes the inner product. 

 Now, consider a connected network with $N$ agents defined by the undirected graph $\mathcal{G}=\{\mathcal{V},\mathcal{E}\}$ with $|\mathcal{V}|=N$ vertices and $|\mathcal{E}|=E$ edges. We assume the graph $\mathcal{G}$ does not contain self-loops and let $L$ denote the graph's Laplacian matrix. Agents try to cooperatively solve the following optimization problem, 
\begin{equation}\label{eq:2.4}
    \min_{x\in \mathbb{C}^K} f(x) := \frac{1}{N}\sum_{i=1}^N f_i(x).
\end{equation}
Let $\mathcal{X}^*$ and $f^*$ denote the optimal set and the minimum function value for optimization Problem \eqref{eq:2.4}. In practice the cross-correlation measurements $d_{ij}^s$ are assumed to be corrupted with additive i.i.d. noise. To recover $x$, each agent $i$ formulates a local least squares error minimization objective function defined as,
\begin{align}\label{eq:2.5}
f_i(x) = \frac{1}{2|\mathcal{N}_i|S_i}\sum_{j\in \mathcal{N}_i} \sum_{s=1}^{S_i} |d_{ij}^s - \langle a_i^s, x \rangle \overline{\langle a_j^s, x \rangle} |^2,
\end{align}
which can be thought of as the sensing process associated with the $i^{th}$ agent of the network~\citep{yazici2019}. Agent $i$'s neighbors are defined as \mbox{$\mathcal{N}_i \triangleq \{ i\in \mathcal{V} | j \in \mathcal{V}, (i,j) \in \mathcal{E} \}$.} The term $S_i$ denotes the number of cross-correlation samples agent $i$ measures. \emph{Thus the graph influences the optimization process and the number of cross-correlation measurements that can be computed.}

Next we state two technical assumptions that will help us later prove convergence for the DGWF algorithm. 
\begin{assum}\label{ass:graph}
The undirected graph $\mathcal{G}$ is connected.
\end{assum}
\begin{assum}\label{ass:set}
The optimal set $\mathcal{X}^*$ of \eqref{eq:2.4}  is nonempty and \mbox{$f^* > -\infty$}. 
\end{assum}

Exploiting the structure of the graph Laplacian $L$, we can thus write the distributed interferometric imaging Problem~\eqref{eq:2.4} equivalently as,
\begin{align}\label{eq:2.7}
    &\min_{\boldsymbol{x}\in \mathbb{C}^{NK}} \Tilde{f}(\boldsymbol{x})\\
    &\textrm{s.t. } \boldsymbol{L}\boldsymbol{x} = \boldsymbol{0}_{NK} \notag,
\end{align}
where $\boldsymbol{x} = \textrm{col}(x_1,...,x_N)$, $\Tilde{f}(\boldsymbol{x}) = \sum_{i=1}^N f_i(x_i)$, and $\boldsymbol{L} = L \otimes I_K$ where $L$ is the graph Laplacian and $\otimes$ is the Kronecker product \citep{uribe2020}. 

\section{Distributed Generalized Wirtinger Flow Algorithm}\label{sec:DGWF}

Recently, \cite{karl2021} presented an algorithm for distributed nonconvex optimization that converges linearly if the local and global cost function is smooth and satisfies the PL inequality, respectively. We propose a primal-dual algorithm inspired by~\cite{karl2021} that uses Wirtinger derivatives ~\citep{yazici2019}. We call the proposed algorithm the \emph{Distributed Generalized Wirtinger Flow} (DGWF) algorithm.

Following the algorithm formulation presented in~\cite{karl2021}, the augmented Lagrangian function for~\eqref{eq:2.7} can be written as,
\begin{equation}\label{eq:3.2}
    \mathcal{A}(\boldsymbol{x},\boldsymbol{u}) = \Tilde{f}(\boldsymbol{x}) +\frac{\lambda_1}{2}\boldsymbol{x}^\top \boldsymbol{L} \boldsymbol{x} + \lambda_2 \boldsymbol{u}^\top \boldsymbol{L}^{1/2} \boldsymbol{x},
\end{equation}
where $\lambda_1 > 0$ and $\lambda_2 > 0$ are the regularization parameters, $\boldsymbol{u}\in \mathbb{C}^{NK}$ is the dual variable, and the constraint $\boldsymbol{L}\boldsymbol{x} = \boldsymbol{0}_{NK}$ is replaced with its equivalent $\boldsymbol{L}^{1/2}\boldsymbol{x} = \boldsymbol{0}_{NK}$ since its nulls are identical~\citep{uribe2020}. \cite{karl2021} proposed the following first order primal-dual gradient method to solve~\eqref{eq:3.2},
\begin{align}
    \boldsymbol{x}_{t+1} &= \boldsymbol{x}_t - \frac{\eta}{\Vert x_0 \Vert} ( \lambda_1 \boldsymbol{L}\boldsymbol{x}_t + \lambda_2 \boldsymbol{L}^{1/2} \boldsymbol{u}_t + \nabla \Tilde{f}(\boldsymbol{x}_t) )\label{eq:3.3.1}\\
    \boldsymbol{u}_{t+1} &= \boldsymbol{u}_t + \frac{\eta}{\Vert x_0 \Vert}\lambda_2 \boldsymbol{L}^{1/2} \boldsymbol{x}_t, \quad \forall \boldsymbol{x}_t \in \mathbb{C}^{NK},\label{eq:3.3.2}
\end{align}
where $\eta > 0$ is the step size and $\Vert x_0 \Vert$ is the norm of the initial estimate. If we denote $\boldsymbol{v}_t = \boldsymbol{L}^{1/2} \boldsymbol{u}_t$, we can rewrite~\eqref{eq:3.3.1} and~\eqref{eq:3.3.2} as,
\begin{align}
    \boldsymbol{x}_{t+1} &= \boldsymbol{x}_t -\frac{\eta}{\Vert x_0 \Vert}(\lambda_1 \boldsymbol{L}\boldsymbol{x}_t + \lambda_2 \boldsymbol{v}_t + \nabla \Tilde{f}(\boldsymbol{x}_t)) \label{eq:3.4.1}\\
    \boldsymbol{v}_{t+1} &= \boldsymbol{v}_t {+} \frac{\eta}{\Vert x_0 \Vert} \lambda_2 \boldsymbol{L} \boldsymbol{x}_t, \forall \boldsymbol{x}_t \in \mathbb{C}^{NK}, \sum_{j=1}^N v_{j,0} = \boldsymbol{0}_{K}. \label{eq:3.4.2}
\end{align}

For initializing the DGWF algorithm we adopt the spectral initialization scheme used in~\cite{yazici2019,yazici2020}. The initial estimate $x_0\in \mathbb{C}^K$ is the rank-1, positive semi-definite matrix approximation of the lifted backprojection estimate \mbox{$\hat{\boldsymbol{X}} \in \mathbb{C}^{K \times K}$},
\begin{align}\label{eq:3.5}
    &\hat{\boldsymbol{X}} {:=} \frac{1}{N}\sum_{i=1}^N \frac{1}{2|\mathcal{N}_i| S_i}\sum_{j\in \mathcal{N}_i} \sum_{s=1}^{S_i} d_{ij}^s a_i^s (a_j^s)^H {+} \overline{d_{ij}^s} a_j^s (a_i^s)^H,
\end{align}
where $x_0 = \sqrt{\lambda_0} \mathbf{v}_0$ with $\lambda_0, \mathbf{v}_0$ being the leading eigenvalue-eigenvector pair of $\hat{\boldsymbol{X}}$ and $(\cdot)^H$ denotes the complex conjugate transpose. The lifted formulation of the $i^{th}$ and $j^{th}$ sensing processes is $F^s = a_j^s (a_i^s)^H $. Let $\boldsymbol{d} = [d_{ij}^1, d_{ij}^2,...,d_{ij}^S] \in \mathbb{C}^S$ denote the vector of cross-correlated measurements. The lifted forward model $\mathcal{F}: \mathbb{C}^{K \times K} \rightarrow \mathbb{C}^S$ is defined as,
\begin{equation}\label{eq:lift}
    \boldsymbol{d} = \mathcal{F}(\hat{\boldsymbol{X}}),
\end{equation}
where $\mathcal{F}$ is a $[S \times K^2]$ matrix with $F^s$ as its rows and $\hat{\boldsymbol{X}}$ concatenated into a vector. This spectral method initializes the iterates within a bounded set i.e., $\mathcal{X}= \{ x \, \vert \, \Vert x \Vert^2 \leq \tau \}$, where $\tau$ is a constant dependent on the restricted isometry constant over rank-1 matrices (RIC$_{\delta_1}$) for the lifted forward model~\eqref{eq:lift}~\cite[Theorem 4.6]{yazici2019}. The initial estimate is distributed to all agents, i.e., $x_0 = x_{1,0} = x_{2,0} =... = x_{i,0}$.
\begin{rem}
The spectral initialization~\eqref{eq:3.5} computes a rank-1 approximation to a low rank matrix recovery problem using all the local cross-correlation measurements stored at every agent. The lifted backprojection estimate $\hat{\boldsymbol{X}}$ is an average over all local lifted backprojection estimates. Thus, one possible way to create a distributed initialization method would be to run a distributed averaging algorithm over the local lifted backprojection estimates \citep{olshevsky2009}. Analysis of how distributed initialization schemes effect the DGWF algorithm is left for future work.
\end{rem}

The DGWF algorithm using spectral initialization \eqref{eq:3.5} and updates \eqref{eq:3.4.1}, \eqref{eq:3.4.2} is shown in pseudo-code as Algorithm~\ref{Algo:1}.
\begin{algorithm}[tb!]
 \caption{Distributed Generalized Wirtinger Flow}
 \begin{algorithmic}[1]
 \renewcommand{\algorithmicrequire}{\textbf{Input:}}
 \renewcommand{\algorithmicensure}{\textbf{Initialize:}}
 \REQUIRE parameters $\lambda_1>0, \lambda_2 >0$, and $\eta>0$
 \ENSURE $x_{i,0} \in \mathbb{C}^K$ with \eqref{eq:3.5} and $v_{i,0} = \boldsymbol{0}_K, \forall i \in [N]$.
  \FOR {$t = 0,1,...$}
  \FOR {$i = 1,...,N$ in parallel}
  \STATE Send $x_{i,t}$ to $\mathcal{N}_i$ and receive $x_{j,t}$ from $j \in \mathcal{N}_i$;
  Primal variable update:
  \STATE $x_{i,t+1} = x_{i,t} - \frac{\eta}{\Vert x_{i,0} \Vert}( \lambda_1 \sum_{j\in \mathcal{N}_i} L_{ij} x_{j,t} + \lambda_2 v_{i,t} + \nabla f_i(x_{i,t}) )$; \\
  Dual variable update:
  \STATE $v_{i,t+1} = v_{i,t} + \frac{\eta}{\Vert x_{i,0} \Vert}\lambda_2 \sum_{j\in \mathcal{N}_i} L_{ij} x_{j,t}$;
  \ENDFOR
  \ENDFOR
 \RETURN $\boldsymbol{x}_t$
 \end{algorithmic}
 \label{Algo:1}
 \end{algorithm}
 
 The gradient $\nabla f_i(x_{i,t})$ is explicitly found by computing the Wirtinger derivative using the local information agent $i$ can access \citep{yazici2019}. For the cross-correlation between agent $i$ and agent $j$ with $(i,j) \in \mathcal{E}$ the gradient is
\begin{align}
    &\nabla f_i(x_{i,t}) {=} \nabla {\left[ \frac{1}{2|\mathcal{N}_i|S_i}\sum_{j\in \mathcal{N}_i}\sum_{s=1}^{S_i} | d_{ij}^s {-} \langle a_i^s,x_{i,t} \rangle \overline{\langle a_j^s,x_{i,t} \rangle}|^2 \right]} \nonumber\\
    &{=} \frac{1}{2|\mathcal{N}_i|S_i}\sum_{j\in \mathcal{N}_i}\sum_{s=1}^{S_i} \overline{e_{ij}^s}\left( a_j^s(a_i^s)^H x_{i,t} \right) {+} e_{ij}^s\left( a_i^s(a_j^s)^H x_{i,t} \right)
\end{align}
where $e_{ij}^s = \left( (a_i^s)^H x_{i,t}(x_{i,t})^H a_j^s - d_{ij}^s \right)$. 

In the next section, we show that when the cost function is smooth and satisfies the PL inequality, the DGWF algorithm converges linearly to a minimum.

\section{Convergence Analysis}\label{sec:conv_analysis}

We follow \cite{yazici2019} geometric analysis of the interferometric imaging problem to show when the lifted forward model's restricted isometry properties over rank-1 positive semi-definite matrices satisfy the Regularity Condition (RC). Then we extend our analysis to show that the RC implies the PL inequality. Finally, convergence of the DGWF algorithm is established.

The RC bounds the local smoothness and curvature of a function, meaning that the function's gradients are \emph{well behaved}.
\begin{defn}\label{def:rc}
(Regularity Condition): A function $f$ satisfies the RC($\alpha$,$\beta$) with $\alpha,\beta > 0$ if,
\begin{equation}\label{eq:condition}
    \langle \nabla f(z),z-x^* \rangle \geq \frac{1}{\alpha} ||\nabla f(z)||^2 + \frac{1}{\beta}||z-x^*||^2,
\end{equation}
for all $z\in \mathbb{R}^{2K}$ and $x^*$ being a minimizer of $f$.
\end{defn}

Next, the PL inequality is defined. The PL inequality does not require convexity but implies that every stationary point is a global minimizer.
\begin{defn}\label{def:pl}
(Polyak-\L{}ojasiewicz Inequality): The function $f$ satisfies the PL inequality if, for some $\mu > 0$,
\begin{equation}
    \frac{1}{2}||\nabla f(z)||^2 \geq \mu (f(z) -f^*), \quad \forall z \in \mathbb{R}^{2K} ,
\end{equation}
where $f:\mathbb{R}^{2K} \rightarrow \mathbb{R}$ is continuously differentiable and $f^*=\min_{z\in \mathbb{R}^{2K}} f(z)$. When this condition holds the function is $\mu$-PL.
\end{defn}

We begin by establishing the relationship between the lifted forward model and the RC using the restricted isometry property (RIP). The RIP arose from the compressed sensing field as a way to measure how close a matrix is to an orthonormal system given a certain degree of sparsity~\citep{candes2005}. 
\begin{defn}\label{def:rip}
(Restricted isometry property): Let $\mathcal{A}:\mathbb{C}^{M \times K} \rightarrow \mathbb{C}^S$ denote a linear operator. Without loss of generality assume $M \leq K$. For every $1 \leq r \leq M$, the $r$-restricted isometry constant (RIC) is defined as the smallest $\delta_r < 1$ such that
\begin{equation}
    (1-\delta_r)\Vert X\Vert_F^2 \leq \Vert \mathcal{A}(X)\Vert^2 \leq (1+\delta_r)\Vert X \Vert_F^2
\end{equation}
holds for all matrices $X$ of rank at most $r$, where \mbox{$\Vert X \Vert_F = \sqrt{Tr(X^H X)}$} denotes the Forbenius norm.
\end{defn}
From~\cite[Theorem 4.6]{yazici2019}, when the lifted forward model $\mathcal{F}$ satisfies the RIP condition over rank-1 positive semi-definite matrices with a RIC$_{\delta_1} < 0.214$ then Definition (\ref{def:rc}) surely holds with any $\alpha, \beta$ satisfying,
\begin{equation}
    \frac{1}{\alpha \Vert x^* \Vert^2} + \frac{c^2(\delta_1)\Vert x^* \Vert^2}{\beta} \leq h(\delta_1) := (1-\delta_2)(1-\epsilon)(2-\epsilon),
\end{equation}
where $\epsilon^2 = (2+\delta_1)(1-\sqrt{1-\frac{\delta_1}{1+\delta_1}}) + \frac{\delta_1^2}{8}$, $\delta_2 = \frac{\sqrt{2}(2+\epsilon)\delta_1}{\sqrt{(1-\epsilon)(2-\epsilon)}}$, and $c(\delta_1) = (2+\epsilon)(1+\epsilon)(1+\delta_1)$.

\subsection{RC \& PL Inequality}\label{sec:rc_pl}
In this subsection we restate~\cite[Lemma 1]{hale2022} in Lemma \ref{lemma:rc_pl} which shows that the RC implies the PL inequality.
\begin{lem}[\cite{hale2022}, Lemma 1]\label{lemma:rc_pl}
Let $f$ \newline have a Lipschitz continuous gradient with Lipschitz constant $L_f$ and the set $\mathcal{X}^* = \{x^*\in \mathbb{R}^{2K} | \nabla f(x^*) = 0\}$ is nonempty and finite. If $f$ is RC($\alpha$,$\beta$), then $f$ is \mbox{${1}/{(\beta^2L_f)}$-PL}.
\end{lem}

\subsection{DGWF Algorithm Convergence Analysis}
Problem \eqref{eq:2.5} is not complex differentiable due to the mapping from $\mathbb{C}^K$ to $\mathbb{R}$~\citep{delgado2009}. In the DGWF algorithm Wirtinger derivatives are used because they provide an elegant way to compute partial derivatives that are differentiable in the complex domain. However, for the convergence analysis it is more convenient to work in the real-valued domain. Thus we equivalently reformulate \eqref{eq:2.5} to the real-valued domain where it is differentiable~\citep{delgado2009}. Without loss of generality, we use~\eqref{eq:2.5} for the $i^{th}$ and $j^{th}$ sensing processes and define $\Tilde{x} \triangleq [\textrm{Re}(x);\textrm{Im}(x)] \in \mathbb{R}^{2K}$, $d_R^s = \textrm{Re}(d_{ij}^s)\in \mathbb{R}$, $d_I^s = \textrm{Im}(d_{ij}^s)\in \mathbb{R}$, and 
\begin{align*}
    \Tilde{A}_R^s \triangleq \begin{bmatrix}
    \Tilde{A}_1^s & \Tilde{A}_2^s \\ 
    -\Tilde{A}_2^s & \Tilde{A}_1^s
    \end{bmatrix} \in \mathbb{R}^{2K\times2K},
    \Tilde{A}_I^s \triangleq \begin{bmatrix}
    \Tilde{A}_2^s & -\Tilde{A}_1^s \\
    \Tilde{A}_1^s & \Tilde{A}_2^s
    \end{bmatrix} \in \mathbb{R}^{2K\times2K},
\end{align*}

with $\Tilde{A}_1^s = \textrm{Re}(a_i^s)^\top \textrm{Re}(a_j^s) + \textrm{Im}(a_i^s)^\top \textrm{Im}(a_j^s)$, and $\Tilde{A}_2^s = -\textrm{Re}(a_i^s)^\top \textrm{Im}(a_j^s) + \textrm{Im}(a_i^s)^\top \textrm{Re}(a_j^s)$, for $s = 1,2,...,S$. We can equivalently rewrite \eqref{eq:2.5} in terms of $\Tilde{x}$ as
\begin{equation}\label{eq:real_fun}
    f(\Tilde{x}) \triangleq \frac{1}{2S}\sum_{s=1}^S \big( \big( d_R^s {-} \Tilde{x}^\top \Tilde{A}_R^s \Tilde{x} \big)^2 {+} \big( d_I^s {-}  \Tilde{x}^\top \Tilde{A}_I^s \Tilde{x} \big)^2 \big).
\end{equation}
\begin{lem} \label{lemma:4}
Function $f(\Tilde{x})$ has Lipschitz gradient on the set $\Tilde{\mathcal{X}} \triangleq \{ \Tilde{x} \, \vert \, \Vert\Tilde{x}\Vert^2 \leq \tau \}$ with Lipschitz constant given by
\begin{align}\label{eq:new_L}
    L_f \triangleq &\frac{1}{S}\sum_{s=1}^S \left( d_R^s \sigma_{\max}(\Tilde{A}_R^s) + d_I^s \sigma_{\max}(\Tilde{A}_I^s) \right. \nonumber \\
    &\left.+ 3\tau \left( \sigma_{\max}^2(\Tilde{A}_R^s) + \sigma_{\max}^2(\Tilde{A}_I^s) \right) \right)  .
\end{align}
\end{lem}

\begin{pf}
Using the definition of the Lipschitz continuity, for $\forall u, v \in \Tilde{\mathcal{X}}$ we have

\begin{align*}
    &\left\Vert \nabla f(u) - \nabla f(v) \right\Vert \\
    &\leq \frac{1}{S}\Bigg\Vert\sum_{s=1}^S \left( (d_R^s - u^\top \Tilde{A}_R^s u)\Tilde{A}_R^s u - (d_R^s - v^\top \Tilde{A}_R^s v)\Tilde{A}_R^s v \right.\\
    &\left.\quad + (d_I^s - u^\top \Tilde{A}_I^s u)\Tilde{A}_I^s u - (d_I^s - v^\top \Tilde{A}_I^s v)\Tilde{A}_I^s v \right) \Bigg\Vert\\
    &\leq \frac{1}{S}\sum_{s=1}^S\left( d_R^s \sigma_{\max}(\Tilde{A}_R^s) + \|u\|^2 \sigma_{\max}^2(\Tilde{A}_R^s) \right) \|u - v\| \\
    &\quad + \frac{1}{S}\sum_{s=1}^S\left\Vert u^\top \Tilde{A}_R^s u\Tilde{A}_R^s v - u^\top \Tilde{A}_R^s v \Tilde{A}_R^s v \right\Vert \\
    &\quad + \frac{1}{S}\sum_{s=1}^S\left\Vert u^\top \Tilde{A}_R^s v\Tilde{A}_R^s v - v^\top \Tilde{A}_R^s v \Tilde{A}_R^s v \right\Vert \\
    &\quad +\frac{1}{S}\sum_{s=1}^S\left( d_I^s \sigma_{\max}(\Tilde{A}_I^s) + \|u\|^2 \sigma_{\max}^2(\Tilde{A}_I^s) \right) \|u - v\| \\
    &\quad + \frac{1}{S}\sum_{s=1}^S\left\Vert u^\top \Tilde{A}_I^s u\Tilde{A}_I^s v - u^\top \Tilde{A}_I^s v \Tilde{A}_I^s v \right\Vert \\
    &\quad + \frac{1}{S}\sum_{s=1}^S\left\Vert u^\top \Tilde{A}_I^s v\Tilde{A}_I^s v - v^\top \Tilde{A}_I^s v \Tilde{A}_I^s v \right\Vert\\
    &\leq \frac{1}{S}\sum_{s=1}^S \left( d_R^s\sigma_{\max}(\Tilde{A_R^s}) + d_I^s\sigma_{\max}(\Tilde{A}_I^s) \right.\\
    &\quad \left.+ 3\tau \left( \sigma_{\max}^2(\Tilde{A}_R^s) + \sigma_{\max}^2(\Tilde{A}_I^s) \right) \right)\Vert u - v\Vert.
\end{align*}
\end{pf}


\begin{prop}
Let the lifted forward model~\eqref{eq:lift} satisfy the RIP condition over
rank-1 positive semi-definite matrices with a RIC$_{\delta_1} \leq 0.214$~\cite[Theorem 4.6]{yazici2019} and Assumptions \ref{ass:graph}-\ref{ass:set} hold. Moreover, let $\{ \tilde{\boldsymbol{x}}_k \}$ be the sequence generated by Algorithm~\ref{Algo:1} with $\lambda_1 > 0$, $\lambda_2 > 0$, $\eta > 0$ applied to~\eqref{eq:real_fun}. Then, there exists a $c_1 > 0$, and $\kappa \in (0,1)$ such that,
\begin{equation}
    \Vert \tilde{\boldsymbol{x}}_t - \bar{\tilde{\boldsymbol{x}}}_t \Vert^2 + N(f(\bar{\tilde{x}}_t) - f^*) \leq (1-\kappa)^t c_1, \quad \forall t \in \mathbb{N}_0.
\end{equation}
where $\bar{\tilde{x}}_t = \frac{1}{N}(\boldsymbol{1}_N^\top \otimes \boldsymbol{I}_K)\boldsymbol{\tilde{x}}_t$ and $\bar{\tilde{\boldsymbol{x}}}_t = \boldsymbol{1}_N \otimes \bar{\tilde{x}}_t$.
\end{prop}

\begin{pf}
It follows from Lemma~\ref{lemma:4} that~\eqref{eq:real_fun} is smooth with Lipschitz constant~\eqref{eq:new_L}. Moreover, if the lifted forward model~\eqref{eq:lift} satisfy the RIP condition over
rank-1 positive semi-definite matrices with a RIC$_{\delta_1} \leq 0.214$~\cite[Theorem 4.6]{yazici2019} it holds that~\eqref{eq:real_fun} is RC, and thus the PL inequality holds from~Lemma~\ref{lemma:rc_pl}. Finally, the desired result follows from~\cite[Theorem 2]{karl2021}.
\end{pf}

\begin{rem}
Due to space limitations the constants $\lambda_1$, $\lambda_2$, $\eta$, $c_1$ and $\kappa$ are not explicitly stated. However they can be directly derived from~\cite[Theorem 1]{karl2021}. 
\end{rem}

\section{Numerical Simulations}\label{sec:sims}

We consider the system model in Fig.~\ref{fig:intro}, where $N$ independent receivers encircle a region of interest to be imaged with $p_i^{rx} \in \mathbb{R}^3$ being the spatial position of the $i^{th}$ receiver. A single transmitter of opportunity at a known spatial position $p^{tx} \in \mathbb{R}^3$ illuminates the imaging scene. The imaging scene domain $D$ can be discretized into $K$ voxels at positions $x\in \mathbb{R}^{K\times 3}$. The reflectivity function $\rho(x)\in\mathbb{C}^K$ for the imaging scene is assumed to be nonfluctuating. Assuming the Born approximation (single bounce) and free-space propagation the frequency domain Green's function solution to the Helmholtz equation for receiver $i$ can be written as~\citep{yazici2019},
\begin{equation}\label{eq:7.1}
    h_i(\omega_s) := J(\omega_s)\int_D e^{-i\omega_s\frac{\phi_i(x)}{c_0}}\alpha_i(x,p^{tx})\rho(x)dx,
\end{equation}
where $\omega_s$ is the $s^{th}$ discretized temporal frequency of the transmitted signal, $J(\omega_s)$ is the transmitted signal power, $c_0$ is the wave velocity, $\alpha_i(x,p^{tx})$ models the propagation attenuation and hardware gains, $\rho(x)$ is the scene reflectivity at position $x \in D$. The $\phi_i(x) = |x-p_i^{rx}| + |x-p^{tx}|$ is the bi-static delay term due to propagation. The radar interferometric imaging model for the $s^{th}$  temporal frequency sample of the $i^{th}$ and $j^{th}$ receivers in the multistatic radar setup can be formulated as,
\begin{multline}\label{eq:7.2}
    d_{ij}(\omega_s) = 
    \sum_{k=1}^K \alpha_i(x_k,p^{tx}) e^{-i\omega_s(|x_k-p_i^{rx}| + |x_k-p^{tx}|)/c_0}\rho_k \times \\ 
    J(\omega_s)\overline{J(\omega_s)}\sum_{k'=1}^K \overline{\alpha_j(x_{k'},p^{tx})}e^{i\omega_s(|x_{k'}-p_j^{rx}| + |x_{k'}-p^{tx}|)/c_0}\overline{\rho_{k'}}.
\end{multline}
From \eqref{eq:7.2} the linear sampling vectors $a_i^s$ and $a_j^s$ can be formulated as,
\begin{align}
    a_i^s &= \left[ J(\omega_s)\alpha_i(x_k,p^{tx}) e^{-i\omega_s(|x_k-p_i^{rx}| + |x_k-p^{tx}|)/c_0} \right]_{k=1}^K \label{eq:7.3}\\
    a_j^s &= \left[ J(\omega_s)\alpha_j(x_k,p^{tx}) e^{-i\omega_s(|x_k-p_j^{rx}| + |x_k-p^{tx}|)/c_0} \right]_{k=1}^K.\label{eq:7.4}
\end{align}
Combing \eqref{eq:7.2} - \eqref{eq:7.4} the radar interferometric imaging model can be vectorized as,
\begin{equation}
    d_{ij}(\omega_s) = \langle a_i^s,\rho \rangle \overline{\langle a_j^s, \rho \rangle}, \quad s=1,...,S; i = 1,...,N, j \neq i
\end{equation}

The simulated true scene reflectivity function and multistatic radar setup is shown in Fig. \ref{fig:intro}. Unless otherwise stated the following assumptions and parameters are assumed: Born approximation, $S=64$ frequency samples, $N = 35$ receivers, $t_{max} = 4000$ for GWF, DGWF max iterations, stepsize $\eta_{t} = \min(1-\exp(-t/\tau_0),0.01)$ where $\tau_0 = 3300$, center frequency $12$ GHz, bandwidth $60$ MHz, image is $12\times 12$ $(K=144)$ with spacing $\Delta = 2.4$ m since Fourier resolution limit $\Delta_{res} = 2.5$ m, graph $\mathcal{G}$ is small-world with a $0.1$ probability of connection, $\lambda_1$ and $\lambda_2 = 1$, signal-to-noise ratio is $50$ dB, and transmitter and receiver gains are $100$ dB.

The mean square error (MSE): $\textrm{MSE} \triangleq \frac{1}{K}\sum_{k=1}^K(\rho_{i,k} - \rho_{k}^*)^2 $ is used to quantitatively determine the accuracy of the reconstructed image $\rho_i$ to the true image $\rho^*$. For the proposed DGWF algorithm the consensus error $\triangleq \sum_{i=1}^N||\rho_i - \frac{1}{N}\sum_{i=1}^N \rho_i||^2$. For each simulation the DGWF algorithm is compared to the centralized GWF algorithm~\citep{yazici2019,yazici2020} using the same cross-correlation measurements as defined by the graph.


\subsection{Effect of Graph Connectivity}
\begin{figure}[tb!]
     \centering
     \begin{subfigure}[c]{1\columnwidth}
         \centering 
         \includegraphics[width=0.7\columnwidth]{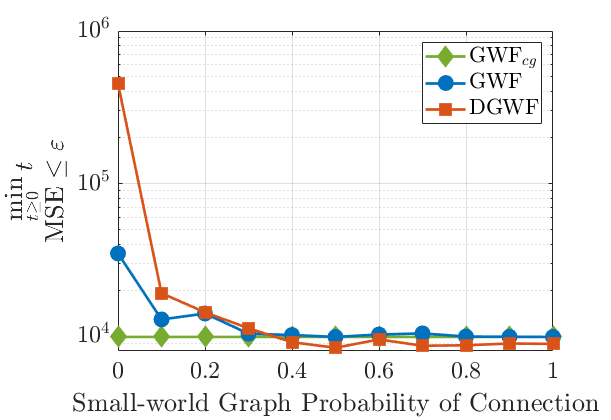}
         \caption{}
         \label{fig:iter_vs_connectivity}
     \end{subfigure}
     \hfill
     \begin{subfigure}[c]{0.46\columnwidth} 
         \centering
         \includegraphics[width=1\columnwidth]{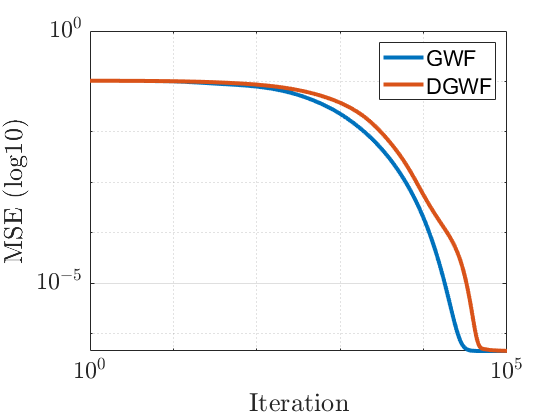}
         \caption{}
         \label{fig:mse_connectivity}
     \end{subfigure}
     \hfill
     \begin{subfigure}[c]{0.46\columnwidth}
         \centering
         \includegraphics[width=1\columnwidth]{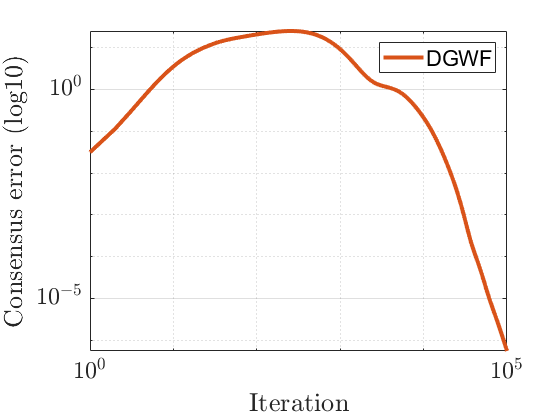}
         \caption{}
         \label{fig:consensus_connectivity}
     \end{subfigure}
        \caption{Comparison between the number of iterations $t$ to reach an MSE $\varepsilon = 10^{-5}$ and the connectivity probability of a small-world graph with $35$ agents (\subref{fig:iter_vs_connectivity}). GWF$_{cg}$ is GWF using a complete graph. MSE (\subref{fig:mse_connectivity}) and consensus error (\subref{fig:consensus_connectivity}) per iteration.}
        \label{fig:effect_graph_connect}
\end{figure}

We examine the effects graph connectivity has on the DGWF algorithm's reconstruction performance. Since the graph specifies how agents communicate and form cross-correlation measurements, we expect the speed of convergence and image reconstruction quality to increase as graph connectivity increases. To test the relationship between graph connectivity and reconstruction performance, we vary the probability of connection for small-world graphs~\citep{strogatz1998} consisting of $35$ agents and record the number of iterations it takes DGWF and GWF to each achieve an MSE of $10^{-5}$. For an additional benchmark, we compare DGWF and GWF against GWF$_{cg}$, which is the GWF algorithm using a complete graph. Figure~\ref{fig:iter_vs_connectivity} shows that decreased graph connectivity results in longer convergence times for GWF and DGWF.
However, after the small-world network's connectivity parameter is greater than $0.4$ there is a negligible change in the number of iterations required to reach the desired reconstruction accuracy.

Additionally, we run the DGWF and GWF algorithms for $10^5$ iterations. Figure~\ref{fig:mse_connectivity} shows that DGWF takes longer to converge compared to GWF because DGWF needs to distribute local information between neighboring agents. The DGWF algorithm reaches a competitive MSE relative to the GWF algorithm without centralized data processing. Figure~\ref{fig:consensus_connectivity} shows that the proposed DGWF algorithm reaches consensus.

\subsection{Effect of Number of Receivers}
\begin{figure}[tb!]
     \centering
     \begin{subfigure}[c]{1\columnwidth}
         \centering 
         \includegraphics[width=0.7\columnwidth]{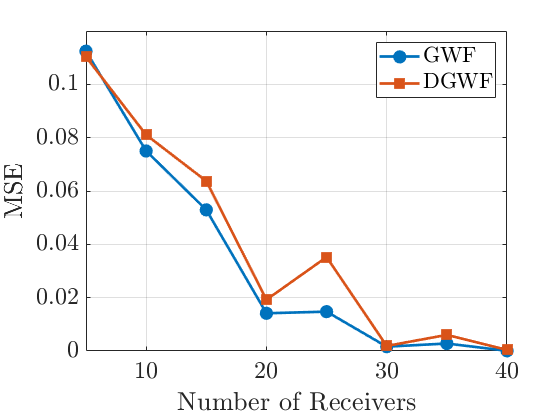}
         \caption{}
         \label{fig:mse_vs_rx}
     \end{subfigure}
     \hfill
     \begin{subfigure}[c]{0.46\columnwidth} 
         \centering
         \includegraphics[width=1\columnwidth]{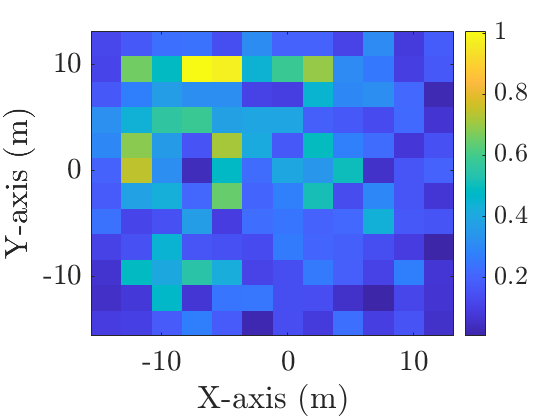}
         \caption{15 receivers.}
         \label{fig:rx_DGWF_15}
     \end{subfigure}
     \hfill
     \begin{subfigure}[c]{0.46\columnwidth}
         \centering
         \includegraphics[width=1\columnwidth]{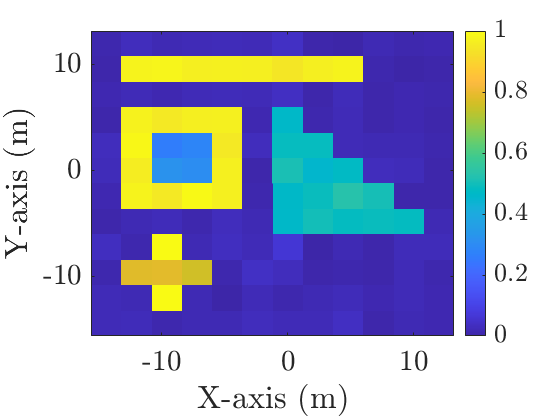}
         \caption{40 receivers.}
         \label{fig:rx_DGWF_40}
     \end{subfigure}
        \caption{Comparison of reconstruction MSE vs. the number of receivers (\subref{fig:mse_vs_rx}). DGWF reconstructions after 4000 iterations for 15 (\subref{fig:rx_DGWF_15}) and 40 (\subref{fig:rx_DGWF_40}) receivers.}
\end{figure}

Next, we analyze the effect the number of receivers has on the convergence behavior for the DGWF and GWF algorithms. We vary the number of agents between 5 and 40 and record the MSE for the DGWF and GWF algorithms after $4000$ iterations. Since the number of receivers directly influences the number of cross-correlation measurements each algorithm can compute, we expect that as the number of receivers decreases, the MSE will increase. In Fig. \ref{fig:mse_vs_rx} the MSE curves for both the DGWF and GWF algorithms show that decreasing the number of receivers increases the reconstruction MSE. Furthermore, Fig. \ref{fig:mse_vs_rx} suggests that after approximately 30 receivers, the addition of more receivers does not significantly decrease the reconstruction error. Example image reconstructions for the DGWF algorithm for 15 and 40 receivers are presented in Fig. \ref{fig:rx_DGWF_15}-\ref{fig:rx_DGWF_40}. It can be observed that when the number of receivers is below 30, the DGWF algorithm reconstructs noticeably degraded images. However, when the number of receivers exceeds 30, the DGWF algorithm achieves significantly accurate image reconstructions.

\section{Conclusion}\label{sec:conclusion}
We studied the distributed interferometric imaging problem over networks. We combine GWF theory~\citep{yazici2019} and distributed optimization methods to create the first distributed interferometric imaging algorithm named the Distributed Generalized Wirtinger Flow (DGWF) algorithm. The DGWF algorithm converges linearly to a global optimum with proper initialization and when the cost function is smooth and satisfies the PL inequality. We demonstrate the capabilities of DGWF through numerical radar interferometric imaging simulations. Future work will study convergence analysis for accelerated methods, time-varying directed graphs, and distributed initialization methods.



\bibliography{ifacconf}             
                                                   







\appendix
\end{document}